\documentclass[12pt]{amsart}

\newcounter{pphcounter}[section]
 \renewcommand{\thepphcounter}{\arabic{pphcounter}}
\newcommand{\pph}{\medskip \refstepcounter{pphcounter}  
 \bf \thepphcounter.~\rm}

\newcommand{\frametitle}[1]{ \pph  \bf{#1}. \rm}
\renewenvironment{frame}{}{}

\usepackage[papersize={210mm, 297mm}]{geometry}
\usepackage{amsmath}

\usepackage[utf8]{inputenc}
\usepackage[english]{babel}

\usepackage{verbatim,latexsym,mathrsfs,color}
\usepackage{tikz}
\usepackage{txfonts}
\usepackage[headings]{fullpage}

\usepackage{comment}

\usepackage{tikz-cd}

 \def\L{L}

\DeclareMathOperator{\Spec}{Spec}

\def\be{\begin{equation}} 
\def\ee{\end{equation}}
\def\bee{\begin{equation}} 
\def\eee{\end{equation}}
\def\bal{\begin{aligned}}   
\def\eal{\end{aligned}}
\def\bes{\begin{equation*}}  
\def\ees{\end{equation*}}
\def\bem{\begin{multline}} 
\def\eem{\end{multline}}

\def\narrower{\advance\leftskip\parindent
  \advance\rightskip\parindent}

\def\Q{\mathbb Q} 
\def\Z{\mathbb Z} 
  
\def\C{\mathbb C} 
  
\def\P{\mathbb P} 
\def\A{\mathbb A}

\def\LL{\mathcal L}

\def\={\;=\;}  
\def\+{\,+\,}

\newcommand{\F}{\mathbb F}

\def\Tr{\operatorname{Tr}}

\def\Spec{\operatorname{Spec}}

\def\ord{\operatorname{ord}}

\def\Bun{\mathrm{Bun}}

\newcommand{\befr}{\begin{frame}}
\newcommand{\efr}{\end{frame}}

\renewcommand{\phi}{\varphi}

\title[Non-abelian Abel's theorems and quaternionic rotation]{Non-abelian Abel's theorems \\ and quaternionic rotation}
\author{V. Golyshev, A. Mellit, V. Rubtsov, and D. van Straten}  

\begin{document}
\bibliographystyle{acm}

\newpage  

\color{black}

\clearpage \setcounter{page}{1}

\begin{flushright}

\end{flushright}

\bigskip
\bigskip
\bigskip

\maketitle

\begin{flushright}
{\it 
}
\end{flushright}
\bigskip
\begin{center}
\parbox{12cm}{\bf Abstract. \rm
\small In order to compute with $l$--adic 
sheaves or crystals 
on a line over~$\F _q$ a low-technology alternative 
to the traditional computation with 
the Hecke operators on the automorphic side
could be helpful. 
A program which has evolved
over the years in our discussions 
with M. Kontsevich 
centers around the concept that, 
in the geometric case, there must exist
certain multiplication laws
on the Galois--representation side
that could be thought of as
precursors of the automorphic lifts: non-abelian 
Abel's theorems, and their restrictions 
to diagonal, Clausen identities. To a varying extent,
they can determine the trace functions of $l$--adic sheaves 
or crystals with prescribed ramification without 
directly appealing to the Hecke--eigen property
on the automorphic side.}
\end{center}

\bigskip
\bigskip

\setcounter{section}{1}

\bigskip

\begin{frame}\frametitle{Motivation} 
The central calculation of this paper 
says that the convolution of the Markov local 
system on the punctured genus $1$ curve 
with the quadratic character sheaf is
a quaternionic rank $3$ local system whose monodromy
is given by the conjugation action
by units in an order in $A_{\{2,3\}}$.
We use it, and related facts,  
as a proof of concept, or illustration, for 
a few ideas 
about Galois representations that are specific 
to the function-field or geometric (= curve over~$\C$)
situation --- namely, that: 
\begin{enumerate}

\item the classical `per prime' formulation of 
lifting by Langlands might be too restrictive, and formulations
seeking to `link' the Euler factors at different 
places should be sought for~\cite{Kontsevich2009}; 

\item Taylor's formula, rather than being an 
isolated fact, is merely the simplest member 
of the family of `master 
formulas' that relate Hecke--type 
operators to hamiltonians
in various quantizations of differential equations;

\item a proper formulation of Abel's theorem 
should take on the multiplicative form, and in such 
form it is capable of surviving in the non-abelian 
situation, potentially furnishing an answer 
to the problem of computing Galois representations;

\item Clausen--type formulas \cite{Katz2009}  should be viewed, 
depending on the optic,  as (precursors of) 
Langlands's lifts or as restrictions to diagonal of 
the multiplication kernels of non-abelian Abel's 
theorems; nevertheless, even these restrictive 
Clausen formulas impose strong conditions 
on the trace functions, and can be used instead 
of the Hecke--eigen 
property on the automorphic side for practically 
computing with them;      

\item in the Betti rendering \cite{BenZviNadler2018}, the presence
of a generalized lift can be signaled by the existence
of a map between character varieties.

\end{enumerate}
The conceptual difference from the geometric Langlands 
correspondence is that we are somewhat agnostic as to 
the nature of the quantization model; non-abelian Abel's theorems
arise when the external powers of the DE are considered.    

We also say a few words about our original motivation 
to classify $D4$ congruence sheaves in the context 
of mirror symmetry. We work up from the

\end{frame}
\begin{frame}

{\underline{Fact}.} 
Classification of Picard rank 1 Fano threefolds 
\cite{IskovskikhProkhorov1999} is
mirrored by the classification of $D3$ 
Kugo-Sato 3-folds  \cite{Golyshev2007}. 
Namely, Fano 3-folds of index $d$ 
and the anticanonical degree $(-K)^3 = 2d^2 N$ correspond $1:1$ to $KS_{N,d}$'s
as in the top right corner of the diagram

\[\begin{tikzcd}
\mathcal{E} \times \mathcal{E}^N \arrow[no head]{r}{} \arrow{d}{} & \widetilde{\mathcal{E} \times \mathcal{E}^N } \arrow{d}{} \\
X_{0}(N)^N \arrow[leftarrow]{r}{\; \; t \mapsto t^d} & \mathbb P^{1} 
\end{tikzcd}
\]

such that the Picard–Fuchs DE 
arising from the right arrow 
is of $D3$ type.
\end{frame}
\begin{frame}

\frametitle{DN equations}
A DN equation \cite{GolyshevStienstra2007} is obtained from  an $(N+1)\times(N+1)$ matrix $A=(a_{ij})_{i,j=0}^N$ 
that satisfies
\be{}\bal
&a_{ij}\=0\,, \quad i-j>1 \\
&a_{ij}\=1\,, \quad i-j=1 \\
&a_{ij}\=a_{N-j,N-i}\,, \quad i-j<1 \\
\eal\ee
The respective differential operator is then defined as 
\[
\L_{A}(t) \= D^{-1} \det\,_{\mathrm{right}} \,\Bigl( \delta_{ij} D  - a_{ij} \bigl( Dt \bigr)^{j-i+1} \Bigr) 
\]
where $\delta_{ij}$ is the Kronecker symbol and $\det_{\mathrm{right}}$ means the right determinant. 
For a sufficiently generic $A$, the respective differential equation 
has maximal unipotent monodromy at $t=0$. The other singularities 
are the inverse eigenvalues of $A$; the respective monodromies 
are orthogonal ($N$ odd) or symplectic ($N$ even) reflections.

\end{frame}
\begin{frame}
Thus, the mirror duals of Picard rank $1$ Fano threefolds are 1-parametric families of motives whose specificity is so strong that it entails modularity.
A literal analogue won't work 
in the non-Shimura situation: 
for $D4$'s, the moduli space of Hodge structures of type 
$$h^{3,0} = h^{2,1} =h^{1,2}=h^{0,3} =1$$
 is 4-dimensional complex, and the universal Hodge structure over it is not a variation of Hodge structures.
\end{frame}
\begin{frame}
We need pencils to play the role of `4-dimensional Kugo-Sato's', but not in any traditional  sense.
Another drawback is that we know essentially nothing about the Calabi-Yau geometries that will appear in the fibers.
\end{frame}
\begin{frame}
A way out is to construct congruence sheaves: representations of $\Q(t)$ with prescribed (small) geometric ramification and congruence properties.
We want to construct them locally over each $\mathbb F_p$, an glue them together across $\Spec \Z$ using congruences as a rigidity constraint.
\end{frame}

\begin{frame}\frametitle{Galois representations and spectral problems}
The prevailing dogma is that in order to compute Galois 
representations, we must be solving certain spectral problems. This is not immediately 
visible in the apparently symmetric Gauss's quadratic reciprocity 
\be
\genfrac(){1pt}{1}{p}{q} = 
\genfrac(){1pt}{1}{q}{p} (-1)^{\frac{p-1}{2} \frac{q-1}{2}},
\ee 
however, the modern interpreation is 
that the meaning of, say, the LHS of 
$\genfrac(){1pt}{1}{5}{q} = \genfrac(){1pt}{1}{q}{5}$ 
is Galois-representational (the action of the $q$-Frobenius 
in $\Q(\!\!\sqrt{5})$), whereas the meaning of the RHS is spectral: 

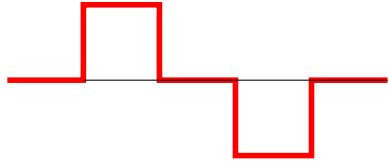
\begin{center}  
\begin{tikzpicture} 
\draw[line width = 0.175em, red] (0,0) -- (1,0) -- (1,1) -- (2,1) -- (2,0) -- (3,0) -- (3,-1) -- (4,-1) -- (4,0) -- (5,0); 
\draw (,0) -- (5,0);

\end{tikzpicture}
\end{center}
define
$
\tilde T_q(f)(x) =\sum_{k=0}^{q-1} f(\frac{x+k}{q}),
$
then the function on the circle shown in the picture is  
$\tilde T_q$--eigen with the eigenvalue $\genfrac(){1pt}{1}{q}{5}$. 
\end{frame}
\begin{frame}
This is clearer in the case of cubic  equations. Put, for instance,  
$f(x)=x^3-4x-1$. There are three ways in which $f$
can split mod prime $p \ne 229 = \mathop{\mathrm{disc}} f$:
with $r_p = 0,1$ or 3 roots in $\F _p$. 
Put 
$a_p = r_p -1, a_{229} = 1$. Extend to $a_n$ by 
multiplicativity in the usual way:
$
\text{ if } \ord_p n = 1, a_n= a_p a_{n/p}$,
$
\text{ if } \ord_p n > 1, a_n= a_p a_{n/p} - \genfrac(){1pt}{1}{229}{q} a_{n/p^2}.
$
\end{frame}
\begin{frame}
Let $K_0(y)$ be the $K$-Bessel function so that
$
K_0(y) = \int_0^\infty \exp (-y \cosh t)\, dt .
$
Finally, define $M$ as the function on the upper half-plane
given by 
$
M(x + iy) = y^{1/2} \sum_{n=1}^\infty a_n \, \cos (nx)\, K_0(ny).
$
\end{frame}
\begin{frame}
Then `reciprocity' means, in particular, that
$
M \left( \frac{i}{229y} \right) = M(iy),
$
which translates into the functional equation
for the respective $L$--function. 
\end{frame}
\begin{frame}

What is the spectral problem in question? Let
$
{\displaystyle \Delta =
-y^{2}\left(\partial ^{2}/\partial x^{2}+
\partial ^{2} / \partial y^{2} \right)}
$
denote the laplacian on the upper half-plane. 
By separation of variables, a series such as above with any 
choice of $a_n$ is $\Delta$-eigen with the eigenvalue $1/4$. 
The deep fact is that $M$ is cuspidal and is $\pm$-invariant under the action of $\Gamma_0(229)$ and, in fact, its normalizer:
$
M\,|_\gamma = \pm M \text{ for } \gamma \in \Gamma_0(229).
$

\end{frame}
\begin{frame}

\end{frame}

\begin{frame}

\bigskip 
The established pattern therefore seems to be: 
in order to compute a Galois representation 
one tries to study a spectral problem 
and recover the Frobenius eigenvalues from the spectral 
parameters. This is true 
without any stretch even historically, e.g.,
\begin{itemize}
\item Ramanujan's $\Delta$ had been found long before the motive 
whose $L$--function is $\Delta$'s Mellin transform \cite{Deligne1973}. 
\item The Poor--Yuen \cite{PoorYuen2015} paramodular form of conductor $61$ over $\Q$  had been found before 
a Calabi--Yau threefold whose $H^3$ gives the respective Galois rep.        

\end{itemize} 

\bigskip

\end{frame}
\begin{frame}

\noindent The spectral problem in question is typically similar to 
the one depicted above:
given a lattice $L$ and a prime number $p$ one considers 
the collection of lattices `mutated at $p$'. For 
$gl(2)$ and rank $2$--lattices these are the lattices
that are of index $p$ in $L$; for a general reductive
group the mutations $m$ correspond to the shape of the
theorem on 
elementary divisors for that group which is in turn
dictated by the structure of its Cartan algebra. 
\end{frame}
\begin{frame}
Averaging over the mutations $m_{s_p}$ of a given shape $s_p$, 
we obtain the action of the Hecke operator $T_{s_p}$ 
on the space of functions on lattices. It turns out that the 
the information contained 
in the collection $\{ T_{s_p} \}$ can be encoded by a semisimple 
conjugacy class in a group $G$ with dual root data (the Satake 
transform), which one can try to interpret as a collection of 
the eigenvalues of the action of $p$-Frobenius in a 
Galois representation.

\end{frame}
\begin{frame}\frametitle{Lifts}
This by now standard worldview suggests that with any group morphism 
$G \to K$ should be associated some transform (lifting) of 
the respective spectral problems for the dual groups.
This was properly codified by Langlands as the 
`functoriality principle' in 1966-7. 
In the number field situation, one way to mediate  
between the automorphic functions
is by considering the $L$--function
\be
L(\rho,s) = \prod_p \frac{1}{\prod_i (1- \alpha^{(p)}_i p^{-s})},
\ee
where $\{ \alpha^{(p)} \}$ is the collection of the eigenvalues
described above. The $L$--functions (as Euler products) 
may undergo simple formal transformations, such as, e.g.
\end{frame}
\begin{frame}
\begin{multline}
\prod_p \frac{1}{(1-\alpha_p p^{-s})(1-p/\alpha_p p^{-s}) }  \rightsquigarrow 
 \prod_p \frac{1}{(1-\alpha_p^2 p^{-s})(1-(p/\alpha_p)^2 p^{-s}) (1-p p^{-s})}.
\end{multline}
for the Shimura $\mathrm{Sym}^2$ lift. It can be expected that these automorphic $L$-
functions should be again automorphic, which can in certain situations be proved with the help of the converse theorems. 
At the level of $L$--functions, such transforms can be represented, roughly, 
by the convolution with theta and Eisenstein kernels. 
There also exist transforms 
at the level of automorphic functions themselves.

\end{frame}

\begin{frame}
It should be noted that $L$--functions lift 
`per individual prime', at the level of each individual 
Euler factor. In other words, one considers 
an essentially non-linear operation at each place, 
and seeks to express it in terms of some integral 
transform/convolution with a kernel. 
There are far fewer examples where one tensors, 
or multiplies, the material at different places. One 
notable example, not of this, but having such flavor, 
is the formula of Gross--Kohnen--Zagier
\cite{GrossKohnenZagier1987} 
that gives integral representation for 
the height pairing of two Heegner points 
corresponding to discriminants $D_1, \, D_2$.      

\end{frame}
\begin{frame} 
A different theta--type transform, 
the Eichler (--Shimizu/Jacquet--Langlands)  
correspondence \cite{Shimizu1977}, establishes a link between Hecke operators acting 
on $2$-lattices of $GL(2)$ and $4$--lattices acted upon by units 
in orders of quaternionic algebras. This is done implicitly 
by using the trace formula; the mechanism in the background is the 
integral transform kernel arising from the Howe duality \cite{Howe1979} between $Gl(2)$
and the quaternions commuting in the metaplectic  representation 
of $GS\! p (4)$. Since the metaplectic 
representation is essentially 
infinite-dimensional and non-geometric, 
it implies that the Eichler correspondence 
should not expected to be given by any geometric kernel. 

\end{frame}
\begin{frame}\frametitle{Lifts as systems of equations on Frobenius traces}
The geometric Langlands correspondence is a `double metaphor'
of this, via the standard chain of abstractions, whereby 
one first passes to the case of function fields over 
the finite field, then to curves over $\C$. The role of
$X_0(229)$ in our example with the Maass form is played by
the space of principal bundles $\Bun_{G^\vee}$; the role of Galois representations 
is played by differential equations. However, the relationship 
between the differential equation on a curve and the 
differential system on $\Bun_{G^\vee}$ is now more direct and has 
geometric nature (is really a correspondence). It follows
that there should exist certains precursors  
of the automorphic lifts at the level of the differential 
equations themselves, which would relate non--linear
operations, such as symmetric/wedge/Schur powers to linear 
operations given by convolutions with kernels.
For $\mathrm{Sym}^2$ we will call these 
`duplication formulas', in contrast to `multiplication 
formulas' where integral representations of 
the result of multiplying solutions 
at different arguments are sought for. (Thus, duplication is similar 
in flavor to Langlands's 
original formulation, while multiplication feels more 
like the Gross-Kohnen-Zagier formula.) 
The fact that the differential equation, or a crystal in the 
$p$-adic setup, has a multiplication kernel ( = satisfies 
a non-abelian Abel's theorem), imposes relations 
on its solutions or the trace functions; 
even duplication alone 
can be a very strong condition. 

\end{frame}
\begin{frame}

\bigskip

The simplest instance of a duplication formula on 
$\A ^1$ reads $(a + x) + x = a + (2 x)$; it can be 
obtained by restricting to diagonal 
the multiplication law  $(a + x) + y = a + (x + y)$, 
which only looks uninspiring if we forget the machinery
of the `master formula' behind it. 

\pph \bf Multiplication laws, Abel's theorems and master formulas. \rm Indeed,  
let $L$ be the differential operator $\frac{d}{dz}$
on $\A ^1$ and consider the spectral problem 
$(L-\lambda) f(z) = 0$. Consider the `normalised' ($f(0)=1$)
solution $\Phi(z, \lambda)$ of this spectral problem, 
$f(z) = \exp (\lambda z) $ and substitute $L$ for
$\lambda$. Then $\exp(\lambda x)\exp (\lambda y) = 
\exp (\lambda (x+y))$ implies 
$\exp(x L)\exp (y L) = 
\exp (\lambda (x+y)L)$, which in turn implies 
an identity between the Hecke/shift operators
$T_xT_y  = T_{x+y}$.

In a higher genus situation, 
multiplication theorems in the abelian setup 
are known as addition theorems, 
but addition is a no-go beyond 
geometric class field theory. 
Abel's discovery 
was that for each algebraic differential $\omega = R(x,y) dx$
there exists a number $P$, such that every sum 
of $N \ge P$  integrals can be reduced to a sum of $P$ integrals.
\bee
 \int_a^{x_1} \omega + \int_a^{x_2} \omega + \dotsb + \int_a^{x_N} \omega 
= \int_a^{y_1} \omega + \int_a^{y_2} \omega + \dotsb + \int_a^{y_P} \omega + E
\eee
where the $ y_1, \dotsc , y_P$ depend  algebraically on
$x_1, \dotsc , x_N$
and $E$ denotes an  elementary function, i.e. rational+$\log$(rational).
\end{frame}
\begin{frame}
An obvious multiplicaticative reformulation 
of Abel's theorem as we know it is 
that 
there exists a simple kernel $K (x | y) $ such that
\bee
\Phi (x_0)\, \Phi (x_1) \dots  \Phi (x_g) = 
\int K(x \, | \, y) \, \Phi (y_1) \dots  \Phi (y_g) \, dy_1 \dots dy_g 
\eee
for any differential equation $d\Phi(x) = \Phi(x) \omega$
on a compact Riemann surface $C$ of genus $g$.
One can fill in all the details; $K$, obviously, is a delta-kernel supported 
on the graph of $\sum (x-o) = \sum (y-o)$  in $\mathrm{Jac}^0(C); \; \Phi (z) = 
\exp \int_o^z \omega$. It is in the multiplicative formulation 
that Abel's theorem is capable of surviving in the non-abelian situation. 

The `master formula' works 
similarly  in the general situation. 
Once a quantization of the DE 
in question is chosen,
the accessory parameters are turned into a commuting
system of hamiltonians on the quantization model. 
The Hecke/shift--type integral operator $T_x$ 
is now obtained by evaluating $\Phi (x)$ with every
accessory parameter replaced by the respective hamiltonian.

\begin{frame}\frametitle{Clausen's formula 
as a duplication formula for Bessel's equation} It runs
\nopagebreak
\be
\left(\sum_{n=0}^\infty \frac{z^n}{n!^2}\right)^2 
= \sum_{n=0}^\infty \frac{z^n}{n!^2}\,\genfrac(){0pt}{0}{2n}{n} = 
1
 + 2\*z
 + \frac{3}{2}\*z^2
 + \frac{5}{9}\*z^3
 + \frac{35}{288}\*z^4
 + \frac{7}{400}\*z^5
+ \cdots 
\ee
and should be interpreted as a $\mathrm{Sym}^2$ lift:  
the LHS is a non-linear transform, whereas 
the RHS is an integral 
representation, namely, by convolution on $\mathbb{G}_m$ 
with a `quadratic character' DE whose solution is 
$
{\displaystyle 
\frac{1}{\sqrt{1-4z}}} 
= \sum_{n=0}^\infty \genfrac(){0pt}{0}{2n}{n}\, z^n
$.  

\bigskip

The relevance of this observation is that 
having a multiplication formula imposes relations on 
trace functions and enables one to recover them. 
To put it even simpler, at the level of the differential equation 
itself, the functional equation    
\be
(\sum_{n=0}^\infty a_n {z^n})^2 
= \sum_{n=0}^\infty \genfrac(){0pt}{0}{2n}{n}\, a_n \, z^n 
\ee
together with the initial condition $a_0 = a_1 =1$ 
enables one to recover the series inductively.  
Bessel's equation is, of course, rigid hypergeometric;
the observation acquires real significance 
in examples with accessory parameters. 
 
\end{frame}

\begin{frame}
\frametitle{Bessel multiplication law} 
The way to pass from Bessel duplication 
to Bessel multiplication is given 
by the Sonine--Gegenbauer formula \cite{GoerlichMarkettStuepp1994},
\cite{Luke1962}.
Define
\be
{\displaystyle J_{0}(z)=
\sum _{m=0}^{\infty }{\frac {(-1)^{m}}{m!^2}}{\left({\frac {z}{2}}\right)}^{2m};}
\ee
it then reads, for $0 < x < y < \infty, \quad y-x  < z < y+x $,
\be
J_0( {x}) \, J_0({y}) = 
\int_{y-x}^{y+x}  K(x,y,z)\,  J_0(z) \, z \, dz,
\ee
with
\bee
K(x,y,z)=\frac{1}{2\pi}\, S(x,y,z)^{-1/2}\eee
where $S(x,y,z)=16^{-1} \, (2x^2y^2+2y^2z^2+2x^2z^2-x^4-y^4-z^4)
$ is the square of the area of a triangle 
of sides $x,y,z$.

\begin{frame}
\frametitle{Bessel's equation as a degenerate case of D2} 
A natural 
idea would be to try to bettify 
 various 
differential equations
and automate the search for lifts for them. 
One expects the lifts for a DE 
to be expressible as identities between the traces 
in the monodromy of the local system of its solutions. 
Since Bessel's equation has an irregular singularity 
at infinity, it would not be convenient to bettify 
the Clausen identity or the Sonine--Gegenbauer
formula directly. Instead, it is useful
to think of Bessel's equation as a degenerate case of 
the generically regular D2 equation. 
Let $
f(t) = t^3+A t^2 + B t,
$
and put 
$
\LL = f(t) \,\partial^2 + f'(t)\, \partial + t, \; 
\partial = \frac{\partial}{\partial t}. 
$
D2~equations are, up to trivial transformation, 
the equations of the form
\bee
\LL\varphi_\lambda(t) = \lambda \varphi_\lambda (t), 
\eee
$\lambda$ being the (only) accessory parameter. 
\end{frame}
\begin{frame}
Define a sequence of polynomials $b_n(\lambda)\in\C[\lambda]$ by
$$
b_0=1, \; \;
b_{n+1}(\lambda)=\frac{1}{B(n+1)^2}\,\left(\bigl(\lambda-A n\, (n+1)\bigr) \, b_n(\lambda) - n^2 b_{n-1}(\lambda)\right).
$$
so that 
$
\varphi_\lambda(t) = \sum_{n=0}^\infty b_n(\lambda) t^n
$
satisfies the D2 equation. Since the degree of $b_n$ is $n$,
it is natural to define 
\end{frame}
\begin{frame}
the numbers $c_{klm}$ by the formula
$
b_k(\lambda) b_l(\lambda) = \sum_{m=0}^\infty c_{klm}\, b_m(\lambda).
$
Expanding, one has
\bee
\sum_{k,l,m=0}^\infty c_{klm} x^k y^l z^m = B\cdot P(x,y,Bz)^{-\frac12},
\eee
where
\bee
P(x,y,z)=(B-xy-yz-xz)^2-4xyz\, (x+y+z+A)
=\mathrm{Discrim}_{\displaystyle{\,t}}\,\Bigl(f(t)-(t-x)\, (t-y) \, (t-z)\Bigr).
\eee
\end{frame}
\begin{frame}
Put
$
K(x,y,z)=B z^{-1} P\left(x,y,B z^{-1}\right)^{-\frac12} \in z^{-1} \C[z^{-1}][[x,y]],
$ cf. \cite{Kontsevich2009}.
One checks that 
$
\LL_x K = \LL_y K = \LL_z K.
$
\end{frame}
\end{frame}
\begin{frame}
For any $\lambda\in\C$ define
\bee
\psi_\lambda (x,y) = \int K(x,y,z) \varphi_\lambda(z) dz.
\eee
We have
$
\LL_x\psi_\lambda = \LL_y\psi_\lambda = \lambda\psi_\lambda,
$
\end{frame}
\begin{frame}
Therefore
$$
\int K(x,y,z) \,\varphi_\lambda(z) \, dz =
 r(\lambda)\, \varphi_\lambda(x)\, \varphi_\lambda(y),
$$
and by substituting $x=0$, from 
$
\varphi_\lambda(0)=1, \, K(0,y,z)=(z-y)^{-1}
$
we obtain \mbox{$r(\lambda)=1$}. 
\newline
Finally, 
the degenerate case $A=B=0$ can be transformed 
into Bessel's equation: \newline
\mbox{$
(\LL-\lambda)\; t^{-1} J_0\, \left(2\sqrt{{-\lambda}/{t}}\right) = 0.
$}
\end{frame}
\begin{frame}
Returning to the Clausen duplication kernel, we see that 
the one arising in the D2 case is
\bee
K(x,x,z) = B\, z^{-1} \Bigl( (x^2-B)^2-4 B \, f(x)\, z^{-1} \Bigr)^{-1/2},
\eee
defining again the convolution with a quadratic character sheaf. 
\medskip
\begin{center}
---
\end{center}
\medskip
\bf Remark. \rm In our language, the differential operator 
$L- \lambda$ is `quantized' tautologically: the model is an 
instance of $\P ^1$, the hamiltonian is $L$ itself, and 
the Hecke--type shift operator at $t_0$ 
is $\varphi_L (t_0)$; the $\Bun_{G^\vee}$--quantization 
coincides with the `separation of variables', or 
Beauville--Mukai, quantization introduced by 
Enriquez and Rubtsov. We refer the reader to 
\cite{EnriquezRubtsov2003}
for the separation-of-variables quantization 
of second--order DEs with more than one accessory parameter. 
        \end{frame}

\end{frame}

\begin{frame}
\frametitle{Bettification and quaternionic rotation} 
Katz's theorem \cite{Katz1996}
says that any two rigid differential 
equations (with r.s.) can be connected by a chain of elementary 
transformations, namely, twists by  character sheaves or 
convolutions with Kummer sheaves w.r. to the additive 
group. 
\end{frame}
\begin{frame}
One wants to expand this type of thinking 
to non--rigid situations. One might attempt to
make the following principle precise 
(to begin with, on a case-by-case basis):

\begin{itemize}

\item Two differential equations with r.s. 
are in the same class if their respective
character varieties are biregular and there is an algebraic transformation formula
between the traces;

\item If two differential equations are 
in the same class, they must be related 
by a chain of convolutions with `elementary' kernels. 
\end{itemize}

\medskip
All rigid DEs 
are in the same class from
this perspective: 
the spaces of accessory parameters are all one-pointed,
and the DEs are indeed related by chains of 
elementary transforms. 
Clausen identities, or `functoriality', 
would become instances of this principle, the 
integral representation for, say, the $\mathrm{Sym}^2$
being one convolution step away from the original~DE. 
However, there will also be more complicated formulas, possibly 
going beyond Langlands's functoriality. 

The main experimental results
of this paper are as follows. Markov's cubic 
\be 
{{\it m_1}}^{2}+{{\it m_2}}^{2}+{{\it m_3}}^{2} = m_1m_2m_3 
\ee 
is a character variety for the free group 
$F_2<A,B>$ realized 
as the fundamental group of a punctured torus
$E \setminus \{O\}$;
we consider those representations 
$\phi: F_2 \longrightarrow SL_2$  
for which the loop around the puncture 
is anti-unipotent:
\bee
\phi ([A,B]) \sim 
\begin{pmatrix} -1 & * \\ 0 & -1 \end{pmatrix}
\eee
and set 
$
m_1 = \Tr \phi(A), \, m_2 = \, \Tr \phi(B),m_3 =\Tr \phi(AB).
$
Let $\mathcal{M}$ be the \emph{Markov local system} on $E$
corresponding to a Markov triple of natural numbers
$(m_1, m_2, m_3)$. 
Consider its `convolution with the quadratic character': 
the fiber at $x$ 
of the local system $\mathcal{M} * \mathcal{L}$ is, by definition, 
$H^1(E,j_*  (\mathcal{M} \otimes \mathcal{L}))$ 
where $\mathcal{L}$ 
is the quadratic character sheaf 
ramified at $2x$ 
and the puncture $O$ and 
uniquely determined by $x$. The local system  
$\mathcal{M} * \mathcal{L}$ is defined on the 
4:1 cover of $E \setminus \{O\}$ that corresponds
to the homomorphism 
$ F_2 \longrightarrow \Z / 2\Z \oplus \Z / 2\Z$
given by $w \mapsto (\mathrm{deg}_Aw , \mathrm{deg}_Bw) \mod 2$. 
Computing, we find experimentally 
:
\begin{enumerate}

\bigskip
\item\sl The traces of the elements $A^2, B^2, (AB)^2$ in the monodromy
representation of 
$\mathcal{M} * \mathcal{L}$ are respectively 
$m_1^2-1,\, m_2^2-1, \, m_3^2-1.$

\medskip

\parshape 1 -0em \dimexpr \linewidth+3em \relax
\noindent  \rm One might therefore think 
that the local system 
$\mathcal{M} * \mathcal{L}$
is, up to a simple
geometric pullback, the $\mathrm{Sym}^2$ of 
$\mathcal{M}$. However,   
 
\smallskip

\parshape 1 3em \dimexpr \linewidth \relax
 \item \sl there are no unipotents in the image of the 
monodromy representation  in~$\mathcal{M} * \mathcal{L}$.

\medskip

\parshape 1 -0em \dimexpr \linewidth+3em \relax
\noindent \rm In fact, 

\smallskip

\parshape 1 3em \dimexpr \linewidth \relax
\item \sl the local system 
$\mathcal{M} * \mathcal{L}$ is a rank 3 local system  
arising from the action by a group of units in 
an order in $A_{\{2,3\}}$ (the quaternion algebra over $\Q$
ramified at $2$ and $3$)
in the conjugation representation 
on traceless quaternions.

\medskip

\parshape 1 -0em \dimexpr \linewidth+3em \relax
\noindent \rm Yet, according to our principle 
stated above, the identity in (1) must not happen   
without there being some longer chain of relations to 
$\mathrm{Sym}^2$, or, a `twisted' Clausen formula of some sort. 
Moreover, 
as integer points are Zariski--dense 
on Markov's surface, one would expect such an identity
to hold for \emph{any} Markov--type local system on $E \setminus
\{O\}$,  i.e. for an arbitrary 
not necessarily integral point 
$(m_1, \, m_2, \, m_3)$.
Indeed,

\smallskip

\parshape 1 3em \dimexpr \linewidth \relax
\item \sl 

let $M$ be \emph{any} Markov--type local system on $E \setminus
\{O\}$, and let $M_{\A^1}$ be the  
Markov-type sheaf on $\A^1$ whose
pullback under the usual double cover
$\sigma : E \setminus
\{O\} \to \A ^1$ is $M$; its monodromies 
around the three finite critical values of 
$\sigma$ are reflections, and the conjugacy 
class of the local monodromy around $\infty$ 
is a size 2 Jordan block with the eigenvalue
$\exp(2 \pi i/4)$.

Let 
$\mathcal{L}_{\chi}$ be 
the Kummer sheaf on $\A ^1$ corresponding 
to a fourth-order character so that 
the local monodromies at $0$ and $\infty$
are $\exp(2 \pi i/4)$ and $\exp(- 2 \pi i/4)$,
and let $\mathcal{L}_{\rho}$ 
be 
the local system on $\A ^1 \setminus 
\{\text{critical values of }\sigma\} $ corresponding 
to an eighth-order character so that 
the local monodromies at the three critical 
values are all equal to $\exp(2 \pi i/8)$. 

One has the following twisted Clausen formula:

$$
\sigma^* \,\Bigl(\mathrm{Sym}^2 \bigl( \mathcal{L}_{\rho} \otimes (M_{\A^1} 
*_{ \mathbb{G}_a} \mathcal{L}_{\chi})\bigr) \Bigr)
\, = \,
\sigma^*(M_{\A^1}) * \mathcal{L} 
$$

\medskip

\parshape 1 -0em \dimexpr \linewidth+3em \relax
\noindent  \rm For the  actual Markov system $\mathcal{M}$, 
the term 
$\mathcal{M}_{\mathbb{A}^1} 
*_{\mathbb{G}_a} \mathcal{L}_{\chi}$
in the 
LHS
can be identified as a subgroup 
of the triangle group (3,4,4), which is known to be
arithmetic, from which fact statement~(3) follows immediately, 
as the RHS is 
$M * \mathcal{L}$.

\bigskip \begin{center}
---
\end{center}
\bigskip 

We thank Maxim Kontsevich for explaining patiently
his `second approach' in \cite{Kontsevich2009}, which 
was the starting point of
our discussions of this subject.

We thank the Institut 
des Hautes \'Etudes Scientifiques in Bures, which has been the venue 
of our meetings on many occasions, for its friendly hospitality. 
VG and VR thank the Institut Henri Poincar\'e for 
support through the Research in Paris program in March 2020.

\end{enumerate}
\end{frame}

\nocite{BeilinsonDrinfeld1996}
\nocite{Langlands1971}
\nocite{Drinfeld1983}
\nocite{Frenkel1995}
\nocite{Frenkel2004}
\nocite{Golyshev2007}
\nocite{GolyshevStienstra2007}
\nocite{GrossZagier1986}
\nocite{GrossKohnenZagier1987}
\nocite{CoatesCortiGalkinGolyshevKasprzyk2013}
\nocite{AlmkvistStratenZudilin2011}
\nocite{IskovskikhProkhorov1999}
\nocite{Katz1996}
\nocite{Katz2009}
\nocite{Kontsevich2009}
\nocite{EtingofKazhdan2021}
\nocite{EtingofFrenkelKazhdan2020}
\nocite{EnriquezRubtsov2003}
\nocite{EnriquezRubtsov2005}
\nocite{Rubtsov2007}
\nocite{BenZviNadler2018}
\nocite{Reiter2015}
\nocite{Shimura1975}
\nocite{Shimizu1977}
\nocite{PoorYuen2015}
\nocite{MillerExton1994}

\bigskip

\bibliographystyle{unsrt}

\bigskip 
\bigskip 

\noindent
{\sc Algebra and Number Theory Laboratory \newline\noindent
Institute for Information Transmission Problems \newline\noindent
Bolshoi Karetny 19, Moscow 127994, Russia}

\nopagebreak 
\smallskip\noindent
golyshev@mccme.ru

\bigskip\noindent
{\sc Faculty of Mathematics
\newline\noindent
University of Vienna
\newline\noindent
Oskar-Morgenstern-Platz 1, Vienna 1090, Austria}

\nopagebreak 
\smallskip\noindent 
anton.mellit@univie.ac.at

\bigskip \noindent
{\sc LAREMA UMR 6093 du CNRS, Mathematics Department, 
\newline\noindent
University of Angers
\newline\noindent
Building I Lavoisier Boulevard Angers, 49045, CEDEX 01, France}

\smallskip\noindent and

\smallskip\noindent
{\sc Theory Division, Institute for Theoretical and Experimental Physics
\newline\noindent
Bolshaya Tcheremushkinskaya, 25, Moscow, 117259, Russia}

\smallskip\noindent
volodya@univ-angers.fr

\bigskip \noindent
{\sc Institut f\"ur Mathematik \newline
Johannes Gutenberg-Universit\"at \newline
Staudingerweg 9, 4. OG
55128 Mainz, Germany
}

\smallskip\noindent
straten@mathematik.uni-mainz.de

\bigskip

\makeatother

\end{document}